\documentclass[a4paper,11pt]{amsart}
\usepackage{fancyhdr}
\usepackage{amsmath}
\usepackage{dsfont}
\usepackage{hyperref}
\usepackage[mathscr]{eucal}
\usepackage[cp1251]{inputenc}
\usepackage[english]{babel}
\usepackage{enumerate,float,indentfirst}
\usepackage{graphicx}
\usepackage{xcolor}
\usepackage{latexsym,a4,mathrsfs,amsthm,amsmath,amssymb,url}
\usepackage{amsfonts}
\usepackage{amssymb}

\usepackage{tikz}
\usepackage{pgfplots}
\usetikzlibrary{plotmarks}

\numberwithin{equation}{section}
\newtheorem{lemma}{Lemma}

\newcommand{\mns}{\scalebox{0.35}[1.0]{\( - \)}}

\begin{document}

\vspace{1in}

\title[Octupoles for octahedral symmetry]{\bf Octupoles for octahedral symmetry}

\author[Yu. Nesterenko]{Yu. Nesterenko}
\address{ Siemens Digital Industries Software }
\email{Yuri.Nesterenko@siemens.com}

\begin{abstract}
Spherical harmonics of degree 4 are widely used in volumetric frame fields design due to their ability to reproduce octahedral symmetry. In this paper we show how to use harmonics of degree 3 (octupoles) for the same purpose, thereby reducing number of parameters and computational complexity. The key ingredients of the presented approach are

\quad \textbullet \ implicit equations for the manifold of octupoles possessing octahedral symmetry up to multiplication by $-1$,

\quad \textbullet \ corresponding rotationally invariant measure of octupole's deviation from the specified symmetry,

\quad \textbullet \ smoothing penalty term compensating the lack of octupoles' symmetries during a field optimization.
\end{abstract}

\maketitle

\thispagestyle{empty}

\section{Introduction}

The most common state-of-the-art approaches for volumetric frame fields design involve the use of 4th-degree spherical harmonics as a representation of field values (see \cite{Huang2011,Ray2015}). This choice is quite natural due to the fact that such harmonics form the linear space closed under 3D rotations and containing harmonics possessing octahedral symmetry. Moreover, algebraic conditions for this symmetry expressed both in terms of implicit equations and penalty function are also known (\cite{Nesterenko2020}).

The major practical disadvantage of this representation is the dimension of the space --- 9, which is three times the minimum required to parameterise frame rotations. In this paper we show how to reduce the dimension by passing from the 4th-degree spherical harmonics to the 3rd-degree ones also known as octupoles (see Landau and Lifshitz \cite{Landau1971} clarifying the connection with the tensor formalism).

In a nutshell, our suggestion is to sacrifice the half of symmetries of frame representation, but compensating for this by additional symmetries in the smoothing penalty term.

\begin{center}
\begin{figure}[h]
\includegraphics[width=13.0cm]{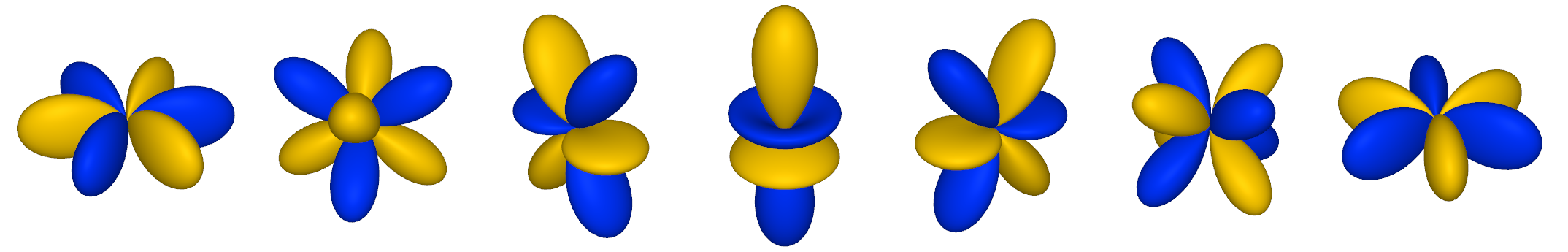}
\caption{ Spherical plots of basis functions $Y_{3,\mns3}, \ldots, Y_{3,3}$. }
\label{fig:basis}
\end{figure}
\end{center}

\section{Semisymmetric octupoles}

As stated, we consider real-valued spherical harmonics of degree 3 on the unit sphere --- octupoles.
These polynomials form the 7D linear space with standard orthonormal basis $Y_{3,\mns3}, \ldots, Y_{3,3}$ (see \cite{Gorller1996}).

With them being odd functions, none of the octupoles (except zero) are octahedrally symmetric by itself, but some of their modules are. Two of such \emph{semisymmetric} octupoles --- $Y_{3,\mns2}$ and $Y_{3,2}$ --- possessing the half of octahedral symmetries (while the remaining are satisfied up to multiplication by $-1$) can be seen in Figure~1.

Since the space we are working in is an eigenspace of the Laplace operator, all possible rotations of its functions lie in the same space and may be represented as a linear combinations of its basis functions.
\begin{center}
\begin{figure}[h!]
\includegraphics[width=8.0cm]{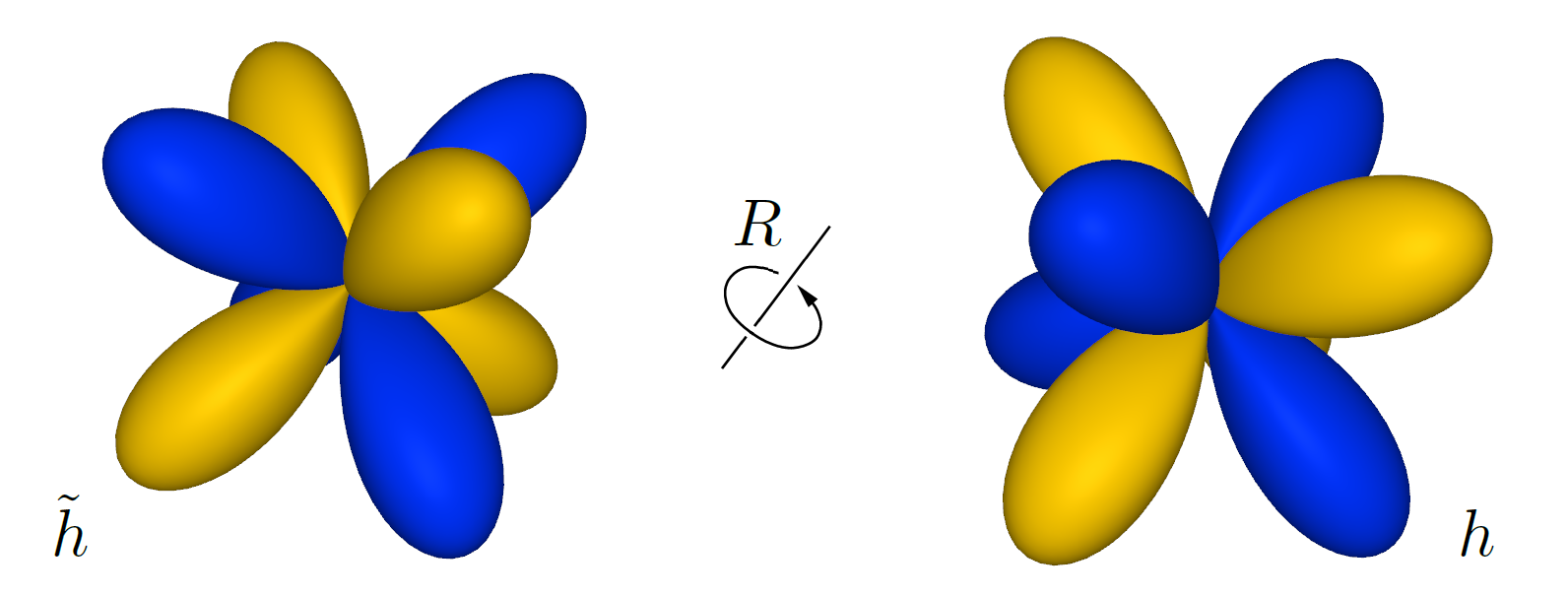}
\caption{ The reference harmonic and its rotation. }
\label{fig:ref7}
\end{figure}
\end{center}
This applies in particular to semisymmetric harmonics. Moreover, in coordinate form all harmonics of this kind may be obtained from a reference one (let it be $Y_{3,\mns2}$) by the formula
\begin{equation*}
a = R_x(\alpha) \times R_y(\beta) \times R_z(\gamma) \times \tilde{a}.
\end{equation*}
Here $\tilde{a} = (0,1,0,0,0,0,0)^T$ and $a$ are coordinates of the reference and the rotated harmonics respectively, and $\alpha$, $\beta$ and $\gamma$ are Euler angles of the corresponding rotation.

Appendix A.1 describes the construction of the rotation matrices $R_x$, $R_y$ and $R_z$.

From the geometrical point of view, all $a(\alpha, \beta, \gamma)$ form
the manifold of dimension 3 embedded in $\mathds{R}^7$. The next lemma claims that this manifold is simply an intersection of quadrics (hypersurfaces of the second order).
\begin{lemma}\label{l1}
The manifold of all semisymmetric octupoles is given by the system of equations
\begin{equation}\label{inp}
\left\{ \begin{split}
&a^T \, a = 1, \\
&a^T \, M_k \, a = 0, \quad k = 1,\ldots,3,
\end{split}\right.
\end{equation}
where $M_1, M_2, M_3$ are the symmetric matrices defined as follows.
\begin{equation}\label{m1}
M_1 = 
\begin{bmatrix}
-5            & 0             & 0             & 0             & 0             & 0             & 0             \\
0             & 0             & 0             & 0             & 0             & 0             & 0             \\
0             & 0             & 3             & 0             & 0             & 0             & 0             \\
0             & 0             & 0             & 4             & 0             & 0             & 0             \\
0             & 0             & 0             & 0             & 3             & 0             & 0             \\
0             & 0             & 0             & 0             & 0             & 0             & 0             \\
0             & 0             & 0             & 0             & 0             & 0             & -5          
\end{bmatrix}
\end{equation}

\begin{equation}\label{m2}
M_2 = 
\begin{bmatrix}
0             & 5             & 0             & 0             & 0             & 0             & 0             \\
5             & 0             & \sqrt{15}     & 0             & 0             & 0             & 0             \\
0             & \sqrt{15}     & 0             & 0             & 0             & 0             & 0             \\
0             & 0             & 0             & 0             & 2             & 0             & 0             \\
0             & 0             & 0             & 2             & 0             & \sqrt{15}     & 0             \\
0             & 0             & 0             & 0             & \sqrt{15}     & 0             & 5             \\
0             & 0             & 0             & 0             & 0             & 5             & 0             
\end{bmatrix}
\end{equation}

\begin{equation}\label{m3}
M_3 = 
\begin{bmatrix}
0             & 0             & 0             & 0             & 0             & 5             & 0             \\
0             & 0             & 0             & 0             & \sqrt{15}     & 0             & -5            \\
0             & 0             & 0             & 2             & 0             & -\sqrt{15}    & 0             \\
0             & 0             & 2             & 0             & 0             & 0             & 0             \\
0             & \sqrt{15}     & 0             & 0             & 0             & 0             & 0             \\
5             & 0             & -\sqrt{15}    & 0             & 0             & 0             & 0             \\
0             & -5            & 0             & 0             & 0             & 0             & 0             
\end{bmatrix}
\end{equation}

\end{lemma}

\emph{Idea of proof.}
The given implicit equations may be obtained by the standard technique based on rational parametrization of the unit circle (to eliminate trigonometric expressions) and Gr\"obner basis construction (see \cite{Cox2015}). 
The statement can be verified by direct calculations.

Topology of the manifold of semisymmetric octupoles is described in Appendix A.2.

\begin{center}
\begin{figure}[h!]
\includegraphics[width=12cm]{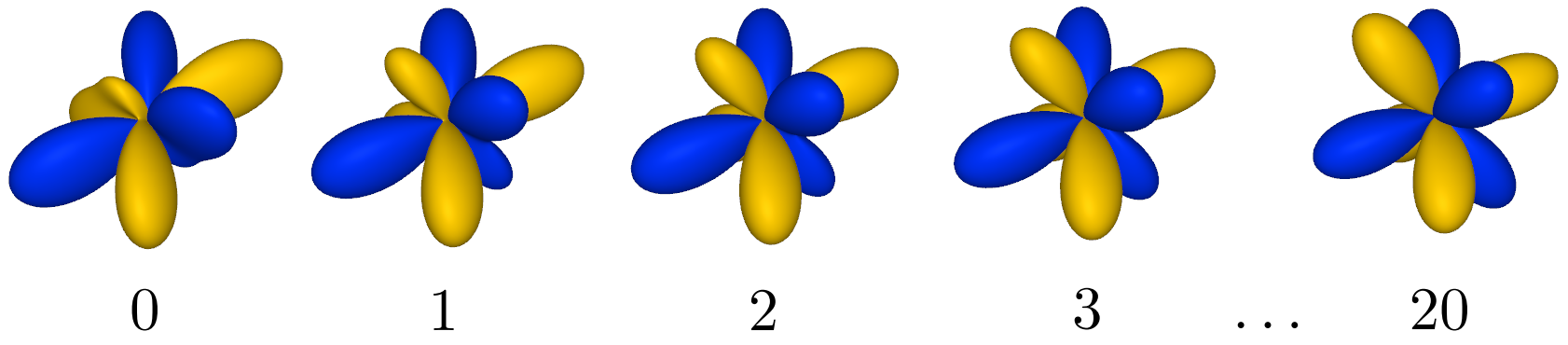}
\caption{ Iterative semisymmetrization }
\label{fig:convergence}
\end{figure}
\end{center}

\section{Semisymmetry enforcement}

Now we are ready to construct the polynomial penalty function enforcing its argument --- octupole --- to be semisymmetric.
\begin{lemma}\label{l2}
Homogeneous 4th-degree polynomial 
\begin{equation}\label{dev}
\begin{split}
d(a) = \,\, & 25 a_{\mns3}^4+50 a_{\mns3}^2 a_{\mns2}^2+20 \sqrt{15} a_{\mns3} a_{\mns2}^2 a_{\mns1}-10 a_{\mns3}^2 a_{\mns1}^2+ \\
       & 30 a_{\mns2}^2 a_{\mns1}^2+ 8 \sqrt{15} a_{\mns3} a_{\mns1}^3+21a_{\mns1}^4-40 a_{\mns3}^2 a_0^2+80 a_{\mns2}^2 a_0^2+ \\
       & 32 a_{\mns1}^2 a_0^2+16 a_0^4+ 120 a_{\mns3} a_{\mns2} a_0 a_1-16 \sqrt{15} a_{\mns2} a_{\mns1} a_0 a_1- \\
       & 10 a_{\mns3}^2a_1^2+30 a_{\mns2}^2 a_1^2- 24 \sqrt{15} a_{\mns3} a_{\mns1} a_1^2+42 a_{\mns1}^2 a_1^2+32 a_0^2 a_1^2+ \\
       & 21 a_1^4+120 a_{\mns3} a_{\mns1} a_0 a_2+ 8 \sqrt{15} a_{\mns1}^2a_0 a_2+40 \sqrt{15} a_{\mns3} a_{\mns2} a_1 a_2- \\
       & 8 \sqrt{15} a_0 a_1^2 a_2+50 a_{\mns3}^2 a_2^2- 20 \sqrt{15} a_{\mns3} a_{\mns1} a_2^2+30 a_{\mns1}^2 a_2^2+ \\
       & 80 a_0^2a_2^2+30 a_1^2 a_2^2-120 a_{\mns2} a_{\mns1} a_0 a_3- 20 \sqrt{15} a_{\mns2}^2 a_1 a_3+ \\
       & 24 \sqrt{15} a_{\mns1}^2 a_1 a_3-8 \sqrt{15} a_1^3 a_3+40 \sqrt{15} a_{\mns2}a_{\mns1} a_2 a_3+ \\
       & 120 a_0 a_1 a_2 a_3+20 \sqrt{15} a_1 a_2^2 a_3+50 a_{\mns3}^2 a_3^2+50 a_{\mns2}^2 a_3^2- \\
       & 10 a_{\mns1}^2 a_3^2-40 a_0^2 a_3^2- 10 a_1^2 a_3^2+50a_2^2 a_3^2+25 a_3^4,
\end{split}
\end{equation}
where $a = (a_{\mns3},a_{\mns2},a_{\mns1},a_{0},a_{1},a_{2},a_{3}) \in \mathds{R}^7$ consists of octupole coordinates in basis $Y_{3,\mns3}, \ldots, Y_{3,3}$, defines the rotationally invariant measure of octupole deviation from semisymmetry.
\end{lemma}

\emph{Idea of proof.}
The statement follows from the method of obtaining this polynomial. It consists of averaging the trial non-invariant deviation measure
\begin{equation}
\widehat{d}(a) = \sum_{k=1}^3 (a^T \, M_k \, a)^2
\end{equation}
over $SO_3$ action's orbits. The next formula can be verified by direct calculations.
\begin{equation}\label{devproof}
\begin{split}
&d(a) \, = \, \frac{1}{\emph{vol} \, SO_3} \int\displaylimits_{R \in SO_3} \widehat{d}(R \cdot a) \, d\mu = \\
 \frac{1}{8\pi^2} \int\displaylimits_0^{2\pi} \int\displaylimits_0^{\pi} & \int\displaylimits_0^{2\pi} \widehat{d}(R_z(\alpha) \times R_x(\beta) \times R_z(\gamma) \times a) \sin\beta \, d\gamma \, d\beta \, d\alpha.
\end{split}
\end{equation}

\section{Numerical example}

In this section we show how deviation measure (\ref{dev}) works.
We use the next combination of the scale and semisymmetry controlling terms with positive weights $w_1 = 5$ and $w_2 = 2.5$
\begin{equation}\label{energy}
p(a; w_1, w_2) = w_1 (a^T \, a - 1)^2 + w_2 \, d(a),
\end{equation}
as the penalty function together with a simple gradient descent method.
Figure \ref{fig:convergence} shows the convergence process of the sample initial octupole to the semisymmetrical one.

The plots below describe the distance to the nearest semisymmetrical octupole and square root of the penalty value. One can see distance-like behavior of the square root of $p(a; w_1, w_2)$.

\quad

\begin{center}
\begin{tikzpicture}
    \begin{axis}[
        height=0.45\textwidth,
        width=0.7\textwidth,
        xlabel=Iterations,
        xtick={0,2,...,20},
        ylabel=Distance measures,
        ymode=log,
       log basis y={10}
]
        \addplot[mark=*,mark options={fill=white},black] coordinates {
(0	,2.81592)
(1	,1.39094)
(2	,1.19559)
(3	,0.985754)
(4	,0.770272)
(5	,0.570106)
(6	,0.403451)
(7	,0.276845)
(8	,0.186271)
(9	,0.123624)
(10	,0.0811893)
(11	,0.0528868)
(12	,0.0342575)
(13	,0.0220825)
(14	,0.0142066)
(15	,0.0091071)
(16	,0.00585892)
(17	,0.00380959)
(18	,0.00242722)
(19	,0.00162)
(20	,0.00096)};

\addplot[mark=triangle*,mark options={fill=white},black] coordinates {
(0	,0.530195)
(1	,0.365971)
(2	,0.303852)
(3	,0.241627)
(4	,0.183293)
(5	,0.132955)
(6	,0.0929796)
(7	,0.0633439)
(8	,0.0423866)
(9	,0.0279998)
(10	,0.0183186)
(11	,0.0119014)
(12	,0.00769476)
(13	,0.00495888)
(14	,0.00318883)
(15	,0.00204835)
(16	,0.0013154)
(17	,0.000845486)
(18	,0.000543288)
(19	,0.000349833)
(20	,2.27E-04)};
    \end{axis}

   	\begin{scope}[shift={(4.5,3.2)}]
	\draw (0,0) --
		plot[mark=triangle*, mark options={fill=white}] (0.25,0) -- (0.5,0)
		node[right]{distance};
	\draw[yshift=\baselineskip] (0,0) --
		plot[mark=*, mark options={fill=white}] (0.25,0) -- (0.5,0)
		node[right]{sqrt penalty};
	\end{scope}
\end{tikzpicture}
\end{center}

Note that due to the invariance of (\ref{energy}) under 3D rotations, its symmetrization effect is orientation agnostic. Therefore, applying it to a field values during the optimization process does not affect octupoles' orientations but helps to maintain their symmetries. This topic is discussed in the next section.

\section{Fields smoothing}

Since octupoles possess octahedral symmetries only up to multiplication by $-1$, we need to compensate for this by defining field smoothness in a special way. For this purpose we propose quite an intuitive expression of the form $|x-y|^2|x+y|^2$ as the smoothing penalty term.

Thus, in discrete cases the final field energy (consideration of boundary conditions is outside the scope of our discussion) becomes
\begin{equation}\label{smoo}
E = \sum_{a \sim b} |a-b|^2|a+b|^2 + \sum_{a} p(a, w_1, w_2).
\end{equation}
Here, the first terms enforce smoothness and the second are responsible for maintenance of the field values semisymmetry during the optimization. 

The described energy function in combination with coarse-to-fine optimization strategy shows results comparable to the "classic" frame fields design approaches (\cite{Huang2011,Ray2015}) while using 7 instead of 9 unknowns per frame.

Two simple examples of this approach are provided in Appendix A.3.

\section{Conclusion}
The implicit equations for the manifold of octupoles possessing octahedral symmetry up to multiplication by $-1$, and the corresponding rotationally invariant deviation measure have been found. The smoothing penalty for octupole fields compensating for the lack of their symmetries has been constructed.

In comparison to existing approaches, the obtained results allow to reduce number of unknowns and computational costs of volumetric frame fields design problems.

\quad

\section{Appendix A.1}
The rotational matrices $R_x$, $R_y$ and $R_z$ for spherical harmonics of degree 3 are defined as follows.

\begin{Small}
\begin{equation}
R_z(\gamma) =
\begin{bmatrix}
\cos 3\gamma  & 0             & 0             & 0             & 0             & 0             &  \sin 3\gamma \\
0             & \cos 2\gamma  & 0             & 0             & 0             &  \sin 2\gamma & 0             \\
0             & 0             & \cos  \gamma  & 0             &  \sin  \gamma & 0             & 0             \\
0             & 0             & 0             & 1             & 0             & 0             & 0             \\
0             & 0             & -\sin  \gamma & 0             & \cos  \gamma  & 0             & 0             \\
0             & -\sin 2\gamma & 0             & 0             & 0             & \cos 2\gamma  & 0             \\
-\sin 3\gamma & 0             & 0             & 0             & 0             & 0             & \cos 3\gamma                
\end{bmatrix}
\end{equation}

\begin{equation}
R_x(\frac{\pi}{2}) = \frac{1}{4}
\begin{bmatrix}
0             & 0             & 0             & \sqrt{10}     & 0             & -\sqrt{6}     & 0             \\
0             & -4            & 0             & 0             & 0             & 0             & 0             \\
0             & 0             & 0             & \sqrt{6}      & 0             & \sqrt{10}     & 0             \\
-\sqrt{10}    & 0             & -\sqrt{6}     & 0             & 0             & 0             & 0             \\
0             & 0             & 0             & 0             & -1            & 0             & -\sqrt{15}    \\
\sqrt{6}      & 0             & -\sqrt{10}    & 0             & 0             & 0             & 0             \\
0             & 0             & 0             & 0             & -\sqrt{15}    & 0             & 1             
\end{bmatrix}
\end{equation}
\end{Small}

\begin{equation}
R_y(\beta) = R_x(\frac{\pi}{2}) \times R_z(\beta) \times R_x(\frac{\pi}{2})^T
\end{equation}

\begin{equation}
R_x(\alpha) = R_y(\frac{\pi}{2})^T \times R_z(\alpha) \times R_y(\frac{\pi}{2})
\end{equation}

See \cite{Blanco1997,Choi1999,Collado1989,Ivanic1996} for more details.

\newpage

\section{Appendix A.2}
The manifold of frame rotations has the topology of the quotient space $SO_3 \, / \, S_4$, where $S_4$ denotes the group of order 24 of all octahedral symmetries. Semisymmetric octupoles are invariant only under even $S_4$ transformations. Hence the corresponding topology is $SO_3 \, / \, A_4$, where $A_4$ denotes the corresponding subgroup.

\begin{center}
\begin{figure}[h]
\begin{tikzpicture}
\label{fig:zones}

\def \d {0.025cm}
\def \r {0.11cm}
\def \rr {0.14cm}

\draw[xshift=160 * \d, yshift=-160 * \d]
[gray!50,dashed]

    (148.5 * \d, 80.3 * \d) -- (0 * \d, -185 * \d) -- (-160 * \d, 57 * \d) -- cycle

;

\draw[xshift=160 * \d, yshift=-160 * \d]
[gray!70,dashed]

    (90 * \d, -60 * \d) -- (130 * \d, -30 * \d)
    (90 * \d, -60 * \d) -- (95 * \d, -95 * \d)
    (90 * \d, -60 * \d) -- (30 * \d, -40 * \d)

    (-10 * \d, 90 * \d) -- (30 * \d, 125 * \d)
    (-10 * \d, 90 * \d) -- (-45 * \d, 120 * \d)
    (-10 * \d, 90 * \d) -- (-10 * \d, 25 * \d)

    (-95 * \d, -70 * \d) -- (-130 * \d, -45 * \d)
    (-95 * \d, -70 * \d) -- (-90 * \d, -105 * \d)
    (-95 * \d, -70 * \d) -- (-45 * \d, -45 * \d)

    (-10 * \d, 25 * \d) -- (30 * \d, -40 * \d)
    (-45 * \d, -45 * \d) -- (-10 * \d, 25 * \d)
    (30 * \d, -40 * \d) -- (-45 * \d, -45 * \d)
;

\begin{scope}[shift={(160 * \d,-160 * \d)}]

\draw[fill=gray!70][gray!70] (0,0) circle (0.33 * \r);
\draw[->,>=stealth][gray!70] (0,0) --  (-60 * \d, -32.5 * \d) node[anchor=south east] {$x$};
\draw[->,>=stealth][gray!70] (0,0) --  (70 * \d, -25 * \d) node[anchor=south west] {$y$};
\draw[->,>=stealth][gray!70] (0,0) --  (-0.1 * \d, 82.5 * \d) node[anchor=west] {$z$};

\end{scope}

\draw[xshift=160 * \d, yshift=-160 * \d][gray!70]

    (-130 * \d, 30 * \d) -- (-130 * \d, -45 * \d) -- (-90 * \d, -105 * \d) -- (-30 * \d, -125 * \d)
    
    (45 * \d, -120 * \d) -- (95 * \d, -95 * \d) --
    (130 * \d, -30 * \d) -- (130 * \d, 45 * \d)
    
    (90 * \d, 105 * \d) -- (30 * \d, 125 * \d) -- (-45 * \d, 120 * \d) -- (-95 * \d, 95 * \d)

    (-90 * \d, 60 * \d) -- (-30 * \d, 40 * \d)

    (10 * \d, -90 * \d) -- (10 * \d, -25 * \d)

    (95 * \d, 70 * \d) -- (45 * \d, 45 * \d)

;

\draw[xshift=160 * \d, yshift=-160 * \d]

    (-95 * \d, 95 * \d) -- (-130 * \d, 30 * \d)

    (-30 * \d, -125 * \d) -- (45 * \d, -120 * \d) 
    
    (130 * \d, 45 * \d) -- (90 * \d, 105 * \d)

    (-90 * \d, 60 * \d) -- (-95 * \d, 95 * \d)

    (-90 * \d, 60 * \d) -- (-130 * \d, 30 * \d)

    (95 * \d, 70 * \d) -- (130 * \d, 45 * \d)
    
    (95 * \d, 70 * \d) -- (90 * \d, 105 * \d)
    
    (10 * \d, -90 * \d) -- (-30 * \d, -125 * \d)
    
    (10 * \d, -90 * \d) -- (45 * \d, -120 * \d)
    
    (10 * \d, -25 * \d) -- (-30 * \d, 40 * \d)
    (45 * \d, 45 * \d) -- (10 * \d, -25 * \d)
    (-30 * \d, 40 * \d) -- (45 * \d, 45 * \d)
;

\begin{scope}[shift={(160 * \d,-160 * \d)}]
\draw[fill=orange] (-130 * \d, 30 * \d) circle (\r);
\draw[gray!70][fill=red!50] (-130 * \d, -45 * \d) circle (\r);
\draw[gray!70][fill=green!50] (-90 * \d, -105 * \d) circle (\r);
\draw[fill=violet!70] (-30 * \d, -125 * \d) circle (\r);
\draw[fill=yellow] (45 * \d, -120 * \d) circle (\r);
\draw[gray!70][fill=orange!50] (95 * \d, -95 * \d) circle (\r);
\draw[gray!70][fill=violet!50] (130 * \d, -30 * \d) circle (\r);
\draw[fill=blue!70] (130 * \d, 45 * \d) circle (\r);
\draw[fill=red] (90 * \d, 105 * \d) circle (\r);
\draw[gray!70][fill=yellow!50] (30 * \d, 125 * \d) circle (\r);
\draw[gray!70][fill=blue!50] (-45 * \d, 120 * \d) circle (\r);
\draw[fill=green] (-95 * \d, 95 * \d) circle (\r);

\draw[fill=violet!70] (-90 * \d, 60 * \d) circle (\r);
\draw[fill=blue!70] (10 * \d, -90 * \d) circle (\r);
\draw[fill=green] (95 * \d, 70 * \d) circle (\r);
\draw[gray!70][fill=green!50] (90 * \d, -60 * \d) circle (\r);
\draw[gray!70][fill=violet!50] (-10 * \d, 90 * \d) circle (\r);
\draw[gray!70][fill=blue!50] (-95 * \d, -70 * \d) circle (\r);

\draw[fill=red] (10 * \d, -25 * \d) circle (\r);
\draw[fill=orange] (45 * \d, 45 * \d) circle (\r);
\draw[fill=yellow] (-30 * \d, 40 * \d) circle (\r);
\draw[gray!70][fill=orange!50] (-10 * \d, 25 * \d) circle (\r);
\draw[gray!70][fill=yellow!50] (-45 * \d, -45 * \d) circle (\r);
\draw[gray!70][fill=red!50] (30 * \d, -40 * \d) circle (\r);

\end{scope}

\draw[xshift=160 * \d, yshift=-160 * \d][thick]

    (0 * \d, 190 * \d) -- (148.5 * \d, 80.3 * \d) -- (160 * \d, -57 * \d) -- (0 * \d, -185 * \d) -- (-148.5 * \d, -80.3 * \d) -- (-160 * \d, 55 * \d) -- cycle

    (0 * \d, 190 * \d) -- (160 * \d, -57 * \d) -- (-148.5 * \d, -80.3 * \d) -- cycle

;

\begin{scope}[shift={(160 * \d,-160 * \d)}]

\draw[fill=gray] (-148.5 * \d, -80.3 * \d) circle (\rr);
\draw[fill=gray] (148.5 * \d, 80.3 * \d) circle (\rr);

\draw[fill=gray] (160 * \d, -57 * \d) circle (\rr);
\draw[fill=gray] (-160 * \d, 57 * \d) circle (\rr);

\draw[fill=gray] (0 * \d, 190 * \d) circle (\rr);
\draw[fill=gray] (0 * \d, -185 * \d) circle (\rr);

\end{scope}

\end{tikzpicture}
\caption{ Fundamental zones for $S_4$ and $A_4$ symmetries. }
\label{fig:rodr}
\end{figure}
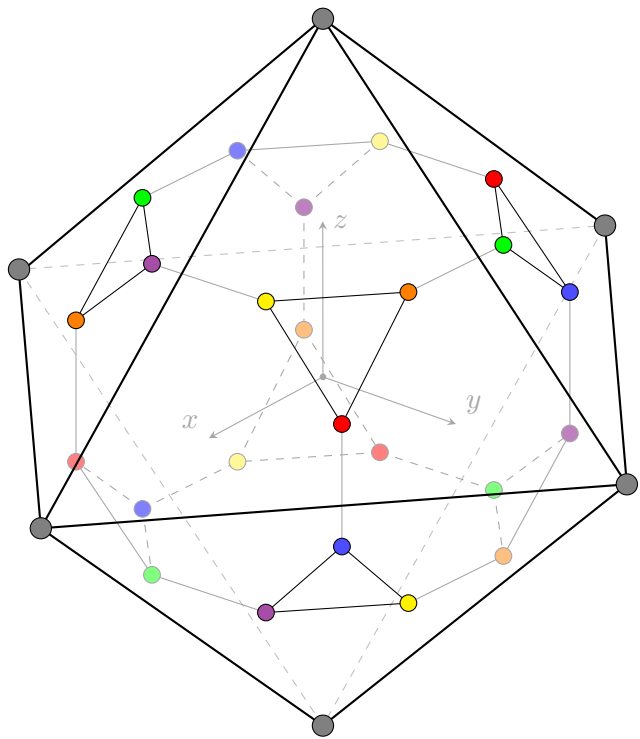
\end{center}

It's possible to describe these topologies using the Rodriguez representation of $3D$ rotations (see \cite{Becker1989}). The fundamental zone for $S_4$ symmetries in this representation has the form of a truncated cube with 6 regular octagonal faces and 8 regular triangular faces. For $A_4$ symmetries the fundamental zone is the regular octahedron. The inclusion $A_4 \subset S_4$ implicates the reverse one for the fundamental zones (see figure \ref{fig:rodr}, note that the triangular faces are pairwise coplanar).

The topologies we consider are obtained by gluing the opposite octagons and the corresponding opposite triangles with $45^\circ$ and $60^\circ$ turn respectively.
The colors in the picture indicate how the vertices map to each other. Note that all octahedron's vertices are coincident and correspond to the octupole opposite to the reference one (i.e. $-Y_{3,\mns2}$).

\newpage

\section{Appendix A.3}
The pictures below show usual singular structures (see \cite{Huang2011,Ray2015}) of frame fields optimized using energy function~(\ref{smoo}).

\begin{center}
\begin{figure}[h!]
\includegraphics[height=5cm]{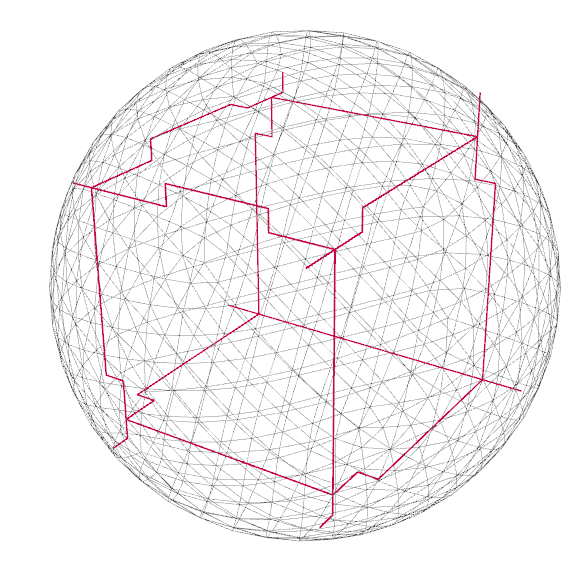}
\includegraphics[height=5cm]{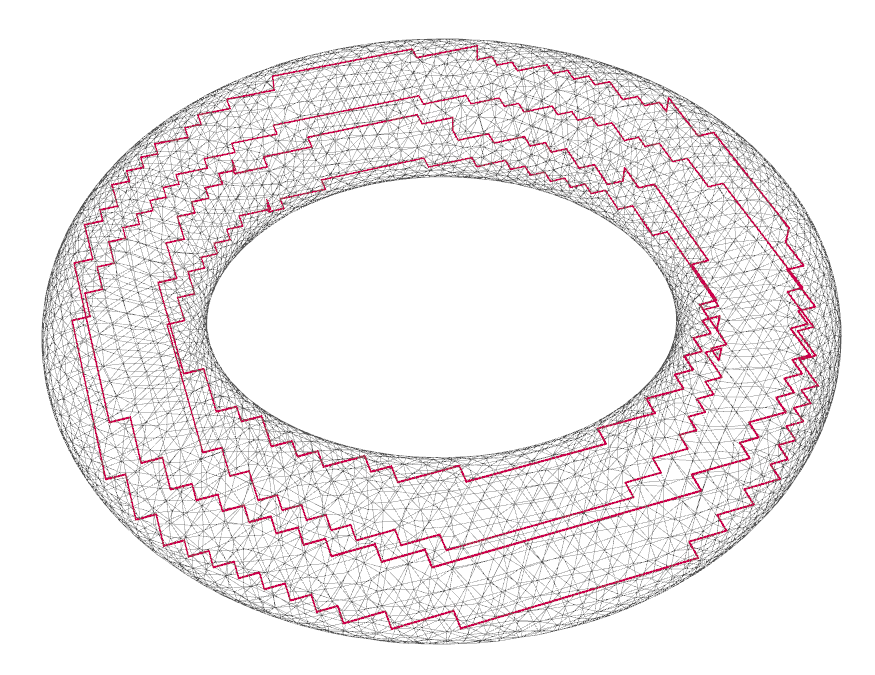}
\caption{ Frame fields singularities of valence 3. }
\label{fig:examples}
\end{figure}
\end{center}

\bibliographystyle{plain}
\bibliography{lit}

\begin{thebibliography}{10}

\bibitem{Becker1989}
R.~Becker and S.~Panchanadeeswaran.
\newblock Crystal rotations represented as {R}odrigues vectors.
\newblock {\em Texture, Stress, and Microstructure}, 10:167--194, 1989.

\bibitem{Blanco1997}
M.A. Blanco, M.~Florez, and M.~Bermejo.
\newblock Evaluation of the rotation matrices in the basis of real spherical
  harmonics.
\newblock {\em Journal of Molecular Structure: THEOCHEM}, 419(1-3):19--27,
  1997.

\bibitem{Choi1999}
C.H. Choi, J.~Ivanic, M.S. Gordon, and K.~Ruedenberg.
\newblock Rapid and stable determination of rotation matrices between spherical
  harmonics by direct recursion.
\newblock {\em The Journal of Chemical Physics}, 111(19):8825--8831, 1999.

\bibitem{Collado1989}
J.R.A. Collado, J.F. Rico, R.~Lopez, M.~Paniagua, and G.~Ramirez.
\newblock Rotation of real spherical harmonics.
\newblock {\em Computer Physics Communications}, 52(3):323--331, 1989.

\bibitem{Cox2015}
D.~O'Shea D.~A.~Cox, J.~Little.
\newblock {\em Ideals, Varieties, and Algorithms}.
\newblock Undergraduate Texts in Mathematics. Springer, forth edition, 2015.

\bibitem{Gorller1996}
C.~G{\"o}rller-Walrand and K.~Binnemans.
\newblock Rationalization of crystal-field parametrization.
\newblock {\em Handbook on the Physics and Chemistry of Rare Earths},
  23:121--283, 1996.

\bibitem{Huang2011}
J.~Huang, Y.~Tong, H.~Wei, and H.~Bao.
\newblock Boundary aligned smooth {3D} cross-frame field.
\newblock {\em ACM Transactions on Graphics}, 30(6):143:1--143:8, 2011.

\bibitem{Ivanic1996}
J.~Ivanic and K.~Ruedenberg.
\newblock Rotation matrices for real spherical harmonics. {D}irect
  determination by recursion.
\newblock {\em The Journal of Chemical Physics}, 100(15):6342--6347, 1996.

\bibitem{Landau1971}
L.~Landau and E.~Lifshitz.
\newblock {\em The Classical Theory of Fields}, pages 96--99.
\newblock Course of Theoretical Physics, Volume 2. Pergamon Press, 3rd revised
  english edition, 1971.

\bibitem{Nesterenko2020}
Yu. Nesterenko.
\newblock On spherical harmonics possessing octahedral symmetry.
\newblock {\em ArXiv}, 2012.12614, 2020.

\bibitem{Ray2015}
N.~Ray and D.~Sokolov.
\newblock On smooth {3D} frame field design.
\newblock {\em ArXiv}, 1507.03351, 2015.

\end{thebibliography}

\end{document}